\documentclass[reqno,12pt]{amsart}

\newtheorem{Theorem}{Theorem}[section]

\newtheorem{Corollary}[Theorem]{Corollary}
\newtheorem{Lemma}[Theorem]{Lemma}
\theoremstyle{remark}

\numberwithin{equation}{section}

\allowdisplaybreaks

\begin{document}

\title[Curious extensions of the ${}_1\psi_1$ sum]
{Curious extensions of Ramanujan's
${}_{\boldsymbol 1}\boldsymbol\psi_{\boldsymbol 1}$ summation formula}

\author[Victor J.\ W.\ Guo]{Victor J.\ W.\ Guo$^*$}
\email{jwguo1977@yahoo.com.cn}
\urladdr{http://math.univ-lyon1.fr/{\textasciitilde}guo}

\author[Michael J.\ Schlosser]{Michael J.\ Schlosser$^{**}$}
\address{Institut f\"ur Mathematik der Universit\"at Wien,
Nordbergstra{\ss}e 15, A-1090 Wien, Austria}
\email{michael.schlosser@univie.ac.at}
\urladdr{http://www.mat.univie.ac.at/{\textasciitilde}schlosse}

\thanks{$^*$The first author was supported by a Junior Research Fellowship
at ESI (Erwin Schr\"odinger Institute), enabling his visit to Vienna
where this research was carried out}
\thanks{$^{**}$The second author was partly supported by
FWF Austrian Science Fund
grants \hbox{P17563-N13}, and S9607 (the second is part
of the Austrian National Research Network
``Analytic Combinatorics and Probabilistic Number Theory'').}
\date{August 1, 2006}

\begin{abstract}
We deduce new $q$-series identities by applying inverse
relations to certain identities for basic hypergeometric series.
The identities obtained themselves do not belong to the hierarchy
of basic hypergeometric series. We extend two of our identities,
by analytic continuation, to bilateral summation formulae which
contain Ramanujan's ${}_1\psi_1$ summation and a
very-well-poised ${}_4\psi_6$ summation as special cases.
\end{abstract}

\maketitle

\section{Introduction}
Ramanujan's ${}_1\psi_1$ summation (cf.~\cite[Eq.~(5.2.1)]{GR})
\begin{equation} \label{11s}
\sum_{k=-\infty}^\infty\frac{(a;q)_k}{(b;q)_k}z^k=
\frac{(q;q)_\infty(b/a;q)_\infty(az;q)_\infty(q/az;q)_\infty}
{(b;q)_\infty(q/a;q)_\infty(z;q)_\infty(b/az;q)_\infty},
\end{equation}
where $|q|<1$ and $|b/a|<|z|<1$, and where
$(x;q)_\infty=\prod_{j\ge 0}(1-xq^j)$
and $(x;q)_k=(x;q)_\infty/(xq^k;q)_\infty$ for integer $k$,
is one of the most important and beautiful identities in
the theory of basic hypergeometric series, see \cite{GR}.

Concerning hypergeometric and basic hypergeometric identities,
there is a dual hierarchy of certain identities closely related to
these but which themselves do not belong to the hierarchy of
hypergeometric or basic hypergeometric series. These identities
can be obtained by applying {\em inverse relations} to the respective
(basic) hypergeometric identities.
For instance, ``dual'' to the {\em binomial theorem},
\begin{equation*}
(a+c)^n=\sum_{k=0}^n\binom nk a^k\,c^{n-k},
\end{equation*}
there is Abel's summation formula,
\begin{equation} \label{abel}
(a+c)^n=\sum_{k=0}^n\binom nk
a(a+bk)^{k-1}(c-bk)^{n-k}
\end{equation}
(cf.~\cite[Sec.~1.5]{R}), which, containing an extra parameter $b$,
is even more general than the binomial theorem.
Similarly, dual to the Chu--Vandermonde summation,
\begin{equation*}
\binom{a+c}n=\sum_{k=0}^n
\binom ak\binom c{n-k},
\end{equation*}
we have the Hagen--Rothe summation,
\begin{equation} \label{rothe}
\binom{a+c}n=\sum_{k=0}^n
\frac{a}{a+bk}\binom{a+bk}k\binom{c-bk}{n-k}
\end{equation}
(cf.~\cite{G}), which, containing an extra parameter $b$, is even more
general than the Chu--Vandermonde summation.
Here we would like to point
out that Abel's summation can be deduced from the Hagen--Rothe
summation. Indeed,
replacing $a$, $b$ and $c$ by $ma$, $mb$ and $mc$, respectively, in
\eqref{rothe},
and dividing both sides by $m^n$, we obtain
\begin{align}
\binom{ma+mc}n m^{-n}
=\sum_{k=0}^n\frac{ma}{ma+mbk}\binom{ma+mbk}k\binom{mc-mbk}{n-k}m^{-n}.
\label{eq:rothe2}
\end{align}
Letting $m\to\infty$ in \eqref{eq:rothe2}, we immediately get
\eqref{abel}.

Furthermore, dual to the Pfaff--Saal\-sch\"utz summation,
\begin{equation}
\frac{(c-a)_n(c-b)_n}{(c)_n(c-a-b)_n}=
\sum_{k=0}^n\frac{(a)_k(b)_k(-n)_k}{(1)_k(c)_k(a+b-c+1-n)_k}
\end{equation}
(cf.~\cite[Thm.~2.2.6]{AAR}), where $(a)_0=1$ and
$(a)_k=a(a+1)\dots(a+k-1)$ for positive integer $k$, we have
the following identity which was derived in \cite[Thm.~7.8]{S1},
\begin{multline}\label{new01gl}
\frac {(2c+1)_n}{(c+1)_n}=
\sum_{k=0}^n\left(\frac {b+(a-c)a} {b+(a-c)(a+k)}\right)
\left(\frac {b+(a+k)^2} {b+a(a+k)}\right)\\\times
\frac {(-n)_k\,(c)_k\,\left(a-c+\frac {b} {a+k}\right)_k}
{(1)_k\,(-c-n)_k\,\left(a+c+\frac {b} {a+k}+1\right)_k}\,
\frac {\left(a+c+\frac {b} {a+k}+1\right)_n}
{\left(a+\frac {b} {a+k}+1\right)_n}.
\end{multline}

Corresponding to the summations in \eqref{abel}, \eqref{rothe}
and \eqref{new01gl}, there exist contiguous identities
(with slightly modified summand, usually involving some additional
linear factors, or ``shifts'' on some of the parameters),
nonterminating summations (expansions),
and basic ($q$-)versions, see \cite{S1}.
For instance, by inverting the $q$-Pfaff--Saalsch\"utz summation,
\begin{equation}
\frac{(c/a;q)_n(c/b;q)_n}{(c;q)_n(c/ab;q)_n}=
\sum_{k=0}^n\frac{(a;q)_k(b;q)_k(q^{-n};q)_k}
{(q;q)_k(c;q)_k(abq^{1-n}/c;q)_k}q^k
\end{equation}
(cf.~\cite[Eq.~(1.7.2)]{GR}),
the following identity was derived in \cite[Thm.~7.34]{S1},
\begin{multline}\label{newq01gl}
\frac {(c^2q;q)_n} {(cq;q)_n}=
\sum_{k=0}^n\frac{(b+(a-c)(a-1))}{(b+(a-c)(a-q^{-k}))}
\frac{(b+(a-q^{-k})^2)}{(b+(a-1)(a-q^{-k}))}\\\times
\frac{(q^{-n};q)_k\,(c;q)_k\,
\big(\frac{b+a(a-q^{-k})}{c(a-q^{-k})};q\big)_k}
{(q;q)_k\,(q^{-n}/c;q)_k\,
\big(cq\frac{b+a(a-q^{-k})}{(a-q^{-k})};q\big)_k}\,
\frac{\big(cq\frac{b+a(a-q^{-k})}{(a-q^{-k})};q\big)_n}
{\big(q\frac{b+a(a-q^{-k})}{(a-q^{-k})};q\big)_n}\,q^k
\end{multline}
(by which we correct some misprints which appeared in the printed
version of \cite{S1}). This identity can be compared to
Equation~\eqref{tsgl} in Theorem~\ref{ts} of this paper,
which is different but somewhat similar to \eqref{newq01gl}. 
To give another example, by inverting the $q$-Gau{\ss} summation,
\begin{equation}
\frac{(c/a;q)_\infty(c/b;q)_\infty}{(c;q)_\infty(c/ab;q)_\infty}=
\sum_{k=0}^\infty\frac{(a;q)_k(b;q)_k}
{(q;q)_k(c;q)_k}\left(\frac c{ab}\right)^k,
\end{equation}
where $|c/ab|<1$ (cf.~\cite[Eq.~(1.5.1)]{GR}), the following
nonterminating identity was derived in \cite[Thm.~7.16]{S1},
\begin{multline}\label{ntnewq0gl}
\frac {(b^2q;q)_{\infty}}{(bq;q)_{\infty}}=
\sum_{k=0}^{\infty}\frac {(c-(a+1)(a+b))}{(c-(a+1)(a+bq^k))}\,
\frac {(c-(a+bq^k)^2)}{(c-(a+b)(a+bq^k))}\\\times
\frac {(b;q)_k\,\big(\frac {(a+bq^k)}{c-a(a+bq^k)};q\big)_k\,
\big(\frac {(a+bq^k)b^2q^{k+1}}{c-a(a+bq^k)};q\big)_{\infty}}
{(q;q)_k\,\big(\frac {(a+bq^k)bq}{c-a(a+bq^k)};q\big)_{\infty}}\,
(bq)^k,
\end{multline}
where $|bq|<1$. This identity can be compared to
Equations~\eqref{tnsgl} and \eqref{tns}, in Theorems~\ref{tns} and
\ref{tnsc}, respectively, of this paper,
which are different but somewhat similar to \eqref{ntnewq0gl}.

In \cite{S1} also {\em multidimensional}\/ identities
associated with root systems of Abel-, Rothe- and the above
``curious'' type are derived. Related beta type integrals
are deduced in \cite{GS} and \cite{S3}.

Concerning {\em bilateral}\/ summations, a ``curious''
generalization of Jacobi's triple product identity was
given in \cite{S2}. However, so far no ``curious'' extensions
of the more general ${}_1\psi_1$ summation formula \eqref{11s}
have been given. Also, so far none of the existing
{\em very-well-poised} summations have been inverted to obtain
non-hypergeometric identities of the above ``curious'' type.
In this paper, we provide for the first time such extensions.
Our formulae, see Theorems~\ref{bns} and \ref{bnsc}, not only
generalize Ramanujan's bilateral summation \eqref{11s}, but also
contain a very-well-poised ${}_4\psi_6$ summation formula as
special case.

Our paper is organized as follows.
In Section~\ref{secpre}, we recall some standard facts about
basic hypergeometric series and list some of the identities
we will be dealing with. In the same section, we also explain
the concept of inverse relations and display some specific
matrix inverses (which are in fact special cases of
Krattenthaler's~\cite{Kr} matrix inverse) we need.
These matrix inverses are utilized in Section~\ref{secqids},
where via inverse relations we deduce from known summations
a couple of new ``curious'' summations which do not belong
to the hierarchy of basic hypergeometric series.
In particular, by inverting the terminating very-well-poised
${}_6\phi_5$ summation we obtain a new terminating ``curious''
summation. Similarly, by inverting the nonterminating
very-well-poised ${}_5\phi_5$ summation we obtain a
new nonterminating ``curious'' summation.
As a limiting case of the new terminating curious
summation, we deduce yet another nonterminating curious
summation. We extend by analytic continuation suitable special
cases of both of these nonterminating unilateral summations
to bilateral summation formulae, which on one hand contain
Ramanujan's ${}_1\psi_1$ summation and on the other hand also
a very-well-poised ${}_4\psi_6$ summation as special cases.

As a matter of fact, we were not able to find a likewise
``curious'' non-hypergeometric generalization of the
${}_6\psi_6$ summation \eqref{66s}, for reasons of convergence.
Such a generalization may still exist but its proof (assuming
it involves inverse relations) would require a matrix inverse
different from Corollaries~\ref{cor1} or \ref{cor2}.
However, already the ``curious'' extensions of Ramanujan's
${}_1\psi_1$ summation in Theorems~\ref{bns} and \ref{bnsc}
came to us as a big surprise.

Concluding this introduction, we would like to add that the
identities derived in this paper have been checked numerically
by {\sc Mathematica}.

\section{Preliminaries}\label{secpre}

\subsection{Basic hypergeometric series}

Let $q$ (the ``base'') be a complex number such that $0<|q|<1$.
Define the {\em $q$-shifted factorial} by
\begin{equation*}
(a;q)_\infty:=\prod_{j\ge 0}(1-aq^j)\qquad\text{and}\qquad
(a;q)_k:=\frac{(a;q)_\infty}{(aq^k;q)_\infty}
\end{equation*}
for integer $k$.
The {\em basic hypergeometric ${}_r\phi_s$ series}
with numerator parameters $a_1,\dots,a_r$, denominator parameters
$b_1,\dots,b_s$, base $q$, and argument $z$ is defined by
\begin{equation*}
{}_r\phi_s\!\left[\begin{matrix}a_1,\dots,a_r\\
b_1,\dots,b_s\end{matrix};q,z\right]:=
\sum_{k=0}^\infty\frac{(a_1;q)_k\dots(a_r;q)_k}
{(q;q)_k(b_1;q)_k\dots(b_s;q)_k}\left((-1)^kq^{\binom k2}\right)^{1+s-r}z^k.
\end{equation*}
The ${}_r\phi_s$ series terminates if one of the numerator parameters
is of the form $q^{-n}$ for a nonnegative integer $n$. If the series
does not terminate, it converges for $r=s+1$ when $|z|<1$.
For $r\le s$, it converges everywhere.
The {\em bilateral basic hypergeometric ${}_r\psi_s$ series}
with numerator parameters $a_1,\dots,a_r$, denominator parameters
$b_1,\dots,b_s$, base $q$, and argument $z$ is defined by
\begin{equation*}
{}_r\psi_s\!\left[\begin{matrix}a_1,\dots,a_r\\
b_1,\dots,b_s\end{matrix};q,z\right]:=
\sum_{k=-\infty}^\infty\frac{(a_1;q)_k\dots(a_r;q)_k}
{(b_1;q)_k\dots(b_s;q)_k}\left((-1)^kq^{\binom k2}\right)^{s-r}z^k.
\end{equation*}
The ${}_r\psi_s$ series reduces to a unilateral
${}_r\phi_s$ series if one of its lower parameters is $q$.
If the series does not terminate, it converges for $r=s$ when $|z|<1$
and $|z|>|b_1\dots b_r/a_1\dots a_r|$.
For $r<s$, it converges everywhere when $|z|>|b_1\dots b_s/a_1\dots a_r|$.

For a thorough exposition on basic hypergeometric series
(or, synonymously, {\em $q$-hypergeometric series}),
including a list of several selected summation and transformation formulas,
we refer the reader to \cite{GR}.

We list some specific identities 
which we will utilize in this paper.

We start with the terminating very-well-poised ${}_6\phi_5$ summation
(cf.\ \cite[Eq.~(2.4.2)]{GR}).
\begin{equation}\label{65s}
{}_6\phi_5\!\left[\begin{matrix}a,\,q\sqrt{a},-q\sqrt{a},b,c,q^{-n}\\
\sqrt{a},-\sqrt{a},aq/b,aq/c,aq^{1+n}\end{matrix}\,;q,
\frac{aq^{1+n}}{bc}\right]=
\frac{(aq;q)_k(aq/bc;q)_k}{(aq/b;q)_k(aq/c;q)_k}.
\end{equation}

This can be extended to the following nonterminating very-well-poised
${}_6\phi_5$ summation (cf.\ \cite[Eq.~(2.7.1)]{GR}):
\begin{multline}\label{65ns}
{}_6\phi_5\!\left[\begin{matrix}a,\,q\sqrt{a},-q\sqrt{a},b,c,d\\
\sqrt{a},-\sqrt{a},aq/b,aq/c,aq/d\end{matrix}\,;q,
\frac{aq}{bcd}\right]\\=
\frac{(aq;q)_\infty(aq/bc;q)_\infty(aq/bd;q)_\infty(aq/cd;q)_\infty}
{(aq/b;q)_\infty(aq/c;q)_\infty(aq/d;q)_\infty(aq/bcd;q)_\infty},
\end{multline}
valid for $|aq/bcd|<1$.
Clearly, \eqref{65ns} reduces to \eqref{65s} for $d=q^{-k}$.

Rather than \eqref{65ns}, we will need the 
following nonterminating very-well-poised
${}_5\phi_5$ summation, resulting from \eqref{65ns} as the special
case where $d\to\infty$ (cf.\ \cite[Ex.~2.22, 2nd Eq.]{GR}):
\begin{equation}\label{55ns}
{}_5\phi_5\!\left[\begin{matrix}a,\,q\sqrt{a},-q\sqrt{a},b,c\\
\sqrt{a},-\sqrt{a},aq/b,aq/c,0\end{matrix}\,;q,
\frac{aq}{bc}\right]=
\frac{(aq;q)_\infty(aq/bc;q)_\infty}
{(aq/b;q)_\infty(aq/c;q)_\infty}.
\end{equation}

The nonterminating ${}_6\phi_5$ summation in \eqref{65ns}
can yet be further extended to Bailey's very-well-poised
${}_6\psi_6$ summation (cf.\ \cite[Eq.~(5.3.1)]{GR}):
\begin{multline}\label{66s}
{}_6\psi_6\!\left[\begin{matrix}q\sqrt{a},-q\sqrt{a},b,c,d,e\\
\sqrt{a},-\sqrt{a},aq/b,aq/c,aq/d,aq/e\end{matrix}\,;q,
\frac{a^2q}{bcde}\right]\\=
\frac{(q;q)_\infty(aq;q)_\infty(q/a;q)_\infty
(aq/bc;q)_\infty}
{(q/b;q)_\infty(q/c;q)_\infty(q/d;q)_\infty(q/e;q)_\infty}\\\times
\frac{(aq/bd;q)_\infty(aq/be;q)_\infty(aq/cd;q)_\infty
(aq/ce;q)_\infty(aq/de;q)_\infty}
{(aq/b;q)_\infty(aq/c;q)_\infty(aq/d;q)_\infty(aq/e;q)_\infty
(a^2q/bcde;q)_\infty},
\end{multline}
valid for $|aq^2/bcde|<1$.
Clearly, \eqref{66s} reduces to \eqref{65ns} for $e=a$.

We will in particular refer to the following
very-well-poised ${}_4\psi_6$ summation formula, obtained
as the $d,e\to\infty$ special case of \eqref{66s}:
 \begin{multline}\label{46s}
{}_4\psi_6\!\left[\begin{matrix}q\sqrt{a},-q\sqrt{a},b,c\\
\sqrt{a},-\sqrt{a},aq/b,aq/c,0,0\end{matrix}\,;q,
\frac{a^2q}{bc}\right]\\=
\frac{(q;q)_\infty(aq;q)_\infty(q/a;q)_\infty
(aq/bc;q)_\infty}
{(q/b;q)_\infty(q/c;q)_\infty(aq/b;q)_\infty(aq/c;q)_\infty}.
\end{multline}

\subsection{Inverse relations}

Let $\mathbb Z$ denote the set of integers
and  $F=(f_{nk})_{n,k\in\mathbb Z}$ be an infinite lower-triangular
matrix; i.e.\ $f_{nk}=0$ unless $n\ge k$.
The matrix $G=(g_{kl})_{k,l\in\mathbb Z}$ is said
to be the {\em inverse matrix} of $F$ if and only if
\begin{equation}\label{orthrel}
\sum_{l\le k\le n} f_{nk}g_{kl}=\delta_{nl}
\end{equation}
for all $n,l\in\mathbb Z$, where $\delta_{nl}$ is the
usual Kronecker delta. Since $F$ anf $G$ are both
lower-triangular, the dual orthogonality relation,
\begin{equation}\label{orthreld}
\sum_{l\le k\le n} g_{nk}f_{kl}=\delta_{nl},
\end{equation}
automatically must hold at the same time.

The method of applying {\em inverse relations} \cite{R}
is a well-known technique for proving identities, or for
producing new ones from given ones.
It is an immediate consequence of the orthogonality relation
\eqref{orthrel}, that if $(f_{nk})_{n,k\in\mathbb Z}$ and
$(g_{kl})_{k,l\in\mathbb Z}$ are lower-triangular matrices
that are inverses of each other, then
\begin{subequations}\label{inv}
\begin{equation}\label{invf}
\sum_{k=0}^nf_{nk}a_k=b_n
\end{equation}
{\em if and only if}
\begin{equation}\label{invg}
\sum_{l=0}^kg_{kl}b_l=a_k.
\end{equation}
\end{subequations}

Another variant of inverse relations which we will also
utilize in this paper involves infinite sums and reads as follows:
\begin{subequations}\label{rotinv}
\begin{equation}\label{rotinvf}
\sum_{n\ge k}f_{nk}a_n=b_k
\end{equation}
{\em if and only if}
\begin{equation}\label{rotinvg}
\sum_{k\ge l}
g_{kl}b_k=a_l,
\end{equation}
\end{subequations}
subject to suitable convergence conditions.
For some applications of \eqref{rotinv} see e.g.\ \cite{Kr,R,S1}.

It is clear that in order to apply \eqref{inv} (or \eqref{rotinv})
effectively, one should have some explicit matrix inversion at hand.

\begin{Lemma}[Krattenthaler~\cite{Kr}]\label{kmi}
Let $(a_j)_{j\in\mathbb Z}$, $(c_j)_{j\in\mathbb Z}$
be arbitrary sequences and $d$ an arbitrary indeterminate. Then
the infinite matrices
$(f_{nk})_{n,k\in\mathbb Z}$ and
$(g_{kl})_{k,l\in\mathbb Z}$ are inverses of each other, where
\begin{subequations}\label{kmigl}
\begin{equation*}
f_{nk}=\frac{\prod_{j=k}^{n-1}(a_j-d/c_k)(a_j-c_k)}
{\prod_{j=k+1}^n(c_j-d/c_k)(c_j-c_k)},
\end{equation*}
\begin{equation*}
g_{kl}=\frac{(a_lc_l-d)(a_l-c_l)}
{(a_kc_k-d)(a_k-c_k)}
\frac{\prod_{j=l+1}^k(a_j-d/c_k)(a_j-c_k)}
{\prod_{j=l}^{k-1}(c_j-d/c_k)(c_j-c_k)}.
\end{equation*}
\end{subequations}
\end{Lemma}

Krattenthaler's matrix inverse is very general as
it contains a vast number of other known explicit infinite matrix
inversions. Several of its useful special cases are of (basic)
hypergeometric type.
The following special case of Lemma~\ref{kmi}
has not been considered explicitly before. It
is exceptional in the sense that although it involves powers of $q$,
it is {\em not} to be considered a $q$-hypergeometric inversion.
(More precisely, the following special case serves as a bridge
between $q$-hypergeometric and certain non-$q$-hypergeometric
identities. For some other such matrix inverses, see \cite{S1}.)

In particular, we set
\begin{equation}
a_j= \frac{1-bc}{1-acq^j},\qquad c_j=1-cq^{-j},\qquad d=1-bc,
\end{equation}
for all integers $j$.

To give a flavor of the elementary computations involved,
we show explicitly how to compute $\prod_{j=l+1}^k(a_j-c_k)$:
\begin{multline*}
\prod_{j=l+1}^k(a_j-c_k)=\prod_{j=l+1}^k\bigg(\frac{1-bc}{1-acq^j}-
\frac {q^k-c}{q^k}\bigg)\\=
\prod_{j=l+1}^k\!\frac{q^k-bcq^k-q^k+c-ac^2q^j+acq^{k+j}}
{(1-acq^j)q^k}=
\prod_{j=l+1}^k\!\frac{c\big(1-bq^k-(c-q^k)aq^j\big)}
{(1-acq^j)q^k}\\=
\left(\frac{c(1-bq^k)}{q^k}\right)^{k-l}
\frac{\Big(\frac{c-q^k}{1-bq^k}aq^{1+l};q\big)_{k-l}}{(acq^{1+l};q)_{k-l}}.
\end{multline*}
Similarly, we compute the other products appearing in \eqref{kmigl}.
After transferring some factors which depend only on one index
from one matrix to the other (which corresponds to simultaneously
multiplying one of the matrices by a suitable diagonal matrix
and multiplying the other matrix by the inverse of that
diagonal matrix) we obtain the following result:

\begin{Corollary}\label{cor1}
Let
\begin{subequations}\label{cor1id}
\begin{multline}
f_{n k}=\frac{(1-bq^n)}{(1-bq^k)}\left(\frac{1-bq^k}{c-q^k}\right)^n
\frac{\Big(1-\frac{c-q^n}{1-bq^n}aq^n\Big)}
{\Big(1-\frac{c-q^k}{1-bq^k}a\Big)}
\frac{\Big(1-\frac{1-bq^k}{c-q^k}q^k\Big)}
{\Big(1-\frac{1-bq^k}{c-q^k}\Big)}\,q^k\\\times
\frac{(q^{-n};q)_k\,(aq^n;q)_k}{(q;q)_k\,(aq;q)_k}
\frac{\Big(\frac{c-q^k}{1-bq^k}a;q\Big)_n}
{\Big(\frac{1-bq^k}{c-q^k}q;q\Big)_n},
\end{multline}
\begin{equation}
g_{k l}=\left(\frac{c-q^k}
{1-bq^k}\right)^l q^{kl}\,
\frac{(1-aq^{2l})}{(1-a)}\frac{(a;q)_l\,(q^{-k};q)_l}{(q;q)_l\,(aq^{1+k};q)_l}
\frac{\Big(\frac{1-bq^k}{c-q^k};q\Big)_l}
{\Big(\frac{c-q^k}{1-bq^k}aq;q\Big)_l}.
\end{equation}
\end{subequations}
Then the infinite matrices $(f_{nk})_{n,k\in\mathbb Z}$ and
$(g_{kl})_{k,l\in\mathbb Z}$ are inverses of each other.
\end{Corollary}

For convenience, we also display another version of Corollary~\ref{cor1},
easily obtained from the above matrix inverse by transferring
some factors from one matrix to the other.

\begin{Corollary}\label{cor2}
Let
\begin{subequations}\label{cor2id}
\begin{equation}
f_{n k}=
\left(\frac{1-bq^k}{c-q^k}q^{-k}\right)^{n-k}
\frac{(1-aq^{2n})}{(1-aq^{2k})}
\frac{(aq^{2k};q)_{n-k}}{(q;q)_{n-k}}
\frac{\Big(\frac{c-q^k}{1-bq^k}aq^k;q\Big)_{n-k}}
{\Big(\frac{1-bq^k}{c-q^k}q^{1+k};q\Big)_{n-k}},
\end{equation}
\begin{multline}
g_{k l}=
(-1)^{k-l}q^{\binom l2-\binom k2}
\frac{(aq;q)_{2k}}{(q;q)_{k-l}\,(aq;q)_{k+l}}\,
\frac{(1-bq^l)}{(1-bq^k)}
\left(\frac{1-bq^k}{c-q^k}\right)^{k-l}\\\times
\frac{\Big(1-\frac{c-q^l}{1-bq^l}aq^l\Big)}
{\Big(1-\frac{c-q^k}{1-bq^k}aq^k\Big)}
\frac{\Big(\frac{c-q^k}{1-bq^k}aq;q\Big)_k}
{\Big(\frac{1-bq^k}{c-q^k};q\Big)_k}
\frac{\Big(\frac{1-bq^k}{c-q^k};q\Big)_l}
{\Big(\frac{c-q^k}{1-bq^k}aq;q\Big)_l}.
\end{multline}
\end{subequations}
Then the infinite matrices $(f_{nk})_{n,k\in\mathbb Z}$ and
$(g_{kl})_{k,l\in\mathbb Z}$ are inverses of each other.
\end{Corollary}

\section{Some curious $q$-series identities}\label{secqids}

We start with a terminating summation, obtained by inverting
the terminating very-well-poised ${}_6\phi_5$ summation
using Corollary~\ref{cor1}.

\begin{Theorem}\label{ts}
Let $a$, $b$, $c$, and $d$ be indeterminates, and let $n$ be
a nonnegative integer. Then 
\begin{multline}\label{tsgl}
\frac{(q/d;q)_n}{(ad;q)_n}\,(ad)^n=
\sum_{k=0}^n\frac{1-bq^n}{1-bq^k}
\left(\frac{1-bq^k}{c-q^k}\right)^n
\frac{(q^{-n};q)_k\,(aq^n;q)_k}{(q;q)_k\,(ad;q)_k}\\\times
\frac{\Big(1-\frac{c-q^n}{1-bq^n}aq^n\Big)}
{\Big(1-\frac{c-q^k}{1-bq^k}aq^k\Big)}
\frac{\Big(1-\frac{1-bq^k}{c-q^k}q^k\Big)}
{\Big(1-\frac{1-bq^k}{c-q^k}\Big)}\,
\frac{\Big(\frac{c-q^k}{1-bq^k}ad;q\Big)_k}
{\Big(\frac{c-q^k}{1-bq^k}a;q\Big)_k}
\frac{\Big(\frac{c-q^k}{1-bq^k}a;q\Big)_n}
{\Big(\frac{1-bq^k}{c-q^k}q;q\Big)_n}\,q^k.
\end{multline}
\end{Theorem}

For $c=1/b$, \eqref{tsgl} reduces to the $q$-Pfaff-Saalsch\"utz summation
(cf.\ \cite[Eq.~(1.7.2)]{GR}). On the other hand, performing
the substitution $c\mapsto -b/c$ and then letting $b\to\infty$
gives the terminating very-well-poised ${}_6\phi_5$ summation \eqref{65s}.

\begin{proof}[Proof of Theorem~\ref{ts}]
Let the inverse matrices $(f_{nk})_{n,k\in\mathbb Z}$ and
$(g_{kl})_{k,l\in\mathbb Z}$ be defined as in Equations~\eqref{cor1id}.
Then \eqref{invg} holds for
\begin{equation*}
a_k=\frac{(aq;q)_k}{(ad;q)_k}
\frac{\Big(\frac{c-q^k}{1-bq^k}ad;q\Big)_k}
{\Big(\frac{c-q^k}{1-bq^k}aq;q\Big)_k}\qquad
\text{and}\qquad
b_l=\frac{(q/d;q)_l}{(ad;q)_l}\,(ad)^l
\end{equation*}
by the terminating ${}_6\phi_5$ summation in \eqref{65s}.
This implies the inverse relation \eqref{invf} with the above values
of $a_k$ and $b_l$. After some minor simplifications we readily
arrive at \eqref{tsgl}.
\end{proof}

Next, we deduce a nonterminating summation, obtained by inverting
the nonterminating very-well-poised ${}_5\phi_5$ summation
using Corollary~\ref{cor2}. 
(We point out that the likewise inversion of the nonterminating
very-well-poised ${}_6\phi_5$ summation fails dues to reasons of
convergence.)

\begin{Theorem}\label{tns}
Let $a$, $b$, $c$, and $d$ be indeterminates. Then 
\begin{multline}\label{tnsgl}
\frac 1{(1-b+a(1-c))}\frac{(ad;q)_\infty}{(aq;q)_\infty}=
\sum_{k=0}^\infty\frac 1{1-bq^k}
\left(\frac{1-bq^k}{c-q^k}\right)^k
\bigg(1-\frac{1-bq^k}{c-q^k}q^k\bigg)\\\times
\frac{(q/d;q)_k}{(q;q)_k\,(aq;q)_k}
\frac{\Big(\frac{c-q^k}{1-bq^k}aq;q\Big)_{k-1}\,
\Big(\frac{1-bq^k}{c-q^k}d;q\Big)_\infty}
{\Big(\frac{1-bq^k}{c-q^k};q\Big)_\infty}\,d^k,
\end{multline}
where $|d/c|<1$.
\end{Theorem}

For $c=1/b$, \eqref{tnsgl} reduces to the $q$-Gau{\ss} summation
(cf.\ \cite[Eq.~(1.5.1)]{GR}). On the other hand, multiplying
both sides by $b$, performing
the substitution $c\mapsto -b/c$ and then letting $b\to\infty$
gives the very-well-poised ${}_5\phi_5$ summation \eqref{55ns}.

\begin{proof}[Proof of Theorem~\ref{tns}]
Let the inverse matrices $(f_{nk})_{n,k\in\mathbb Z}$ and
$(g_{kl})_{k,l\in\mathbb Z}$ be defined as in Equations~\eqref{cor2id}.
Then \eqref{rotinvf} holds for
\begin{equation*}
a_n=(-1)^nq^{\binom n2}
\frac{(q/d;q)_n}{(ad;q)_n}\,d^n
\end{equation*}
and
\begin{equation*}
b_k=(-1)^kq^{\binom k2}d^k\,
(q/d;q)_k
\frac{(aq^{1+2k};q)_\infty}{(ad;q)_\infty}
\frac{\Big(\frac{1-bq^k}{c-q^k}d;q\Big)_\infty}
{\Big(\frac{1-bq^k}{c-q^k}q^{1+k};q\Big)_\infty}
\end{equation*}
by the nonterminating ${}_5\phi_5$ summation in \eqref{55ns}.
This implies the inverse relation \eqref{rotinvg} with the above values
of $a_k$ and $b_l$. After the substitutions $a\mapsto aq^{-2l}$,
$b\mapsto bq^{-l}$, $c\mapsto cq^l$, and $d\mapsto dq^l$, we can get
rid of $l$, and after some simplifications we readily
arrive at \eqref{tnsgl}.
\end{proof}

We now derive a summation ``contiguous'' to Theorem~\ref{tns}.
For this, we simultaneously perform the substitutions $a\mapsto q^{-n}/d$
and $d\mapsto adq^{1+n}$ in \eqref{tnsgl}, multiply both sides by
$(-1)^nq^{\binom n2-n}d^n/(c-d)$ and let $n\to\infty$ (while appealing to
Tannery's theorem~\cite{Br} for justification of taking term-wise limits).
We have the following:

\begin{Theorem}\label{tnsc}
Let $a$, $b$, $c$, and $d$ be indeterminates. Then 
\begin{multline}\label{tnscgl}
\frac 1{(c-d)}\frac{(adq;q)_\infty}{(aq;q)_\infty}=
\sum_{k=0}^\infty\frac 1{c-q^k}
\left(\frac{1-bq^k}{c-q^k}\right)^k
\bigg(1-\frac{1-bq^k}{c-q^k}q^k\bigg)\\\times
\frac{(1/d;q)_k}{(q;q)_k\,(aq;q)_k}
\frac{\Big(\frac{c-q^k}{1-bq^k}aq;q\Big)_k\,
\Big(\frac{1-bq^k}{c-q^k}dq;q\Big)_\infty}
{\Big(\frac{1-bq^k}{c-q^k};q\Big)_\infty}\,d^k,
\end{multline}
where $|d/c|<1$.
\end{Theorem}

Again, for $c=1/b$, \eqref{tnscgl} reduces to the $q$-Gau{\ss} summation.
Also, multiplying both sides by $c$, performing
the substitution $c\mapsto -b/c$ and then letting $b\to\infty$
gives the very-well-poised ${}_5\phi_5$ summation \eqref{55ns}.

Now we extend the $a=0$ case of the unilateral summation in
Theorem~\ref{tns} to a bilateral identity by analytic continuation.
(The likewise extension of the full  summation in
Theorem~\ref{tns} to a bilateral identity by analytic continuation
fails due to reasons of convergence and analyticity.)

\begin{Theorem}\label{bns}
Let $b$, $c$, $d$, and $e$ be indeterminates. Then 
\begin{multline}\label{bnsgl}
\frac 1{(e-b)}
\frac{(q;q)_\infty(de;q)_\infty}{(d;q)_\infty(eq;q)_\infty}
=\sum_{k=-\infty}^\infty\frac 1{(1-bq^k)}
\bigg(1-\frac{1-bq^k}{c-q^k}q^k\bigg)\\\times
\frac{\Big(\frac{1-bq^k}{c-q^k}d;q\Big)_\infty
\Big(\frac{c-q^k}{1-bq^k}eq;q\Big)_\infty}
{\Big(\frac{1-bq^k}{c-q^k};q\Big)_\infty
\Big(\frac{c-q^k}{1-bq^k}q;q\Big)_\infty}\,
\frac{(q/d;q)_k}{(eq;q)_k}
\left(\frac{1-bq^k}{c-q^k}d\right)^k,
\end{multline}
where $|d/c|<1$ and $|e/b|<1$.
\end{Theorem}

Clearly, Theorem~\ref{bns} reduces to the $a=0$ case of
Theorem~\ref{tns} for $e=1$. The $c=1/b$ case of \eqref{bnsgl} gives,
after some rewriting, Ramanujan's ${}_1\psi_1$ summation \eqref{11s}.
On the other hand, performing the substitution $c\mapsto -b/c$ and
then letting $b\to\infty$ gives the very-well-poised
${}_4\psi_6$ summation in \eqref{46s}. We find this unification
of \eqref{11s} and \eqref{46s} quite surprising due to the fact
that the ${}_1\psi_1$ summation is {\em not} a special case
of the very-well-poised ${}_4\psi_6$ summation, in fact, not even
of Bailey's very-well-poised ${}_6\psi_6$ summation formula.

\begin{proof}[Proof of Theorem~\ref{bns}]
We apply Ismail's argument~\cite{Is} to the parameter $e$ using the
$a=0$ case of the nonterminating identity in Theorem~\ref{tns}.
Both sides of \eqref{bnsgl} are analytic in the parameter $e$
in a domain around the origin.
This follows from expanding the products depending on $k$ by
(special cases of) the $q$-binomial theorem.
Now, the identity is true for $e=q^l$, by the
$a=0$ case of Theorem~\ref{tns} (see the next paragraph for the details).
This holds for all $l\ge 0$.
Since $\lim_{l\to\infty}q^l=0$ is an interior point
in the domain of analyticity of $e$, by the identity theorem,
we obtain an identity for general $e$.

The details are displayed as follows. Setting $e=q^l$, the
right-hand side of \eqref{bnsgl} becomes
\begin{multline}\label{longl1}
\sum_{k=-l}^\infty\frac 1{(1-bq^k)}
\bigg(1-\frac{1-bq^k}{c-q^k}q^k\bigg)\\\times
\frac{\Big(\frac{1-bq^k}{c-q^k}d;q\Big)_\infty}
{\Big(\frac{1-bq^k}{c-q^k};q\Big)_\infty
\Big(\frac{c-q^k}{1-bq^k}q;q\Big)_l}\,
\frac{(q/d;q)_k}{(q^{1+l};q)_k}
\left(\frac{1-bq^k}{c-q^k}d\right)^k.
\end{multline}
We shift the summation index in \eqref{longl1} by
$k\mapsto k-l$, and obtain
\begin{multline*}
\frac{(q/d;q)_{-l}}{(q^{1+l};q)_{-l}}\,
\sum_{k=0}^\infty\frac 1{(1-bq^{k-l})}
\bigg(1-\frac{1-bq^{k-l}}{c-q^{k-l}}q^{k-l}\bigg)\\\times
\frac{\Big(\frac{1-bq^{k-l}}{c-q^{k-l}}d;q\Big)_\infty}
{\Big(\frac{1-bq^{k-l}}{c-q^{k-l}};q\Big)_\infty
\Big(\frac{c-q^{k-l}}{1-bq^{k-l}}q;q\Big)_l}\,
\frac{(q^{1-l}/d;q)_k}{(q;q)_k}
\left(\frac{1-bq^{k-l}}{c-q^{k-l}}d\right)^{k-l}\\=
\frac{(q;q)_l}{(d;q)_l}\,q^{-l}\,
\sum_{k=0}^\infty\frac 1{(1-bq^{k-l})}
\bigg(1-\frac{1-bq^{k-l}}{c-q^{k-l}}q^{k-l}\bigg)\\\times
\frac{\Big(\frac{1-bq^{k-l}}{c-q^{k-l}}d;q\Big)_\infty}
{\Big(\frac{1-bq^{k-l}}{c-q^{k-l}}q^{-l};q\Big)_\infty}\,
\frac{(q^{1-l}/d;q)_k}{(q;q)_k}
\left(\frac{1-bq^{k-l}}{c-q^{k-l}}d\right)^k.
\end{multline*}
Next, we apply the $a\mapsto 0$, $b\mapsto bq^{-l}$, $c\mapsto cq^l$,
$d\mapsto dq^l$, case of Theorem~\ref{tns}, simplify, and obtain
for the last expression
\begin{equation*}
\frac{(q;q)_l}{(d;q)_l}\frac{q^{-l}}{(1-bq^{-l})}=
\frac 1{(q^l-b)}\frac{(q;q)_\infty(dq^l;q)_\infty}
{(d;q)_\infty(q^{1+l};q)_\infty},
\end{equation*}
which is exactly the $e=q^l$, case of the
left-hand side of \eqref{bnsgl}.
\end{proof}

Similarly, we can extend the $a=0$ case of the unilateral summation in
Theorem~\ref{tnsc} to a bilateral identity by analytic continuation.
(Again, the likewise extension of the full  summation in
Theorem~\ref{tnsc} to a bilateral identity by analytic continuation
fails due to reasons of convergence and analyticity.)

\begin{Theorem}\label{bnsc}
Let $b$, $c$, $d$, and $e$ be indeterminates. Then 
\begin{multline}\label{bnscgl}
\frac 1{(c-d)}
\frac{(q;q)_\infty(deq;q)_\infty}{(dq;q)_\infty(eq;q)_\infty}
=\sum_{k=-\infty}^\infty\frac 1{(c-q^k)}
\bigg(1-\frac{1-bq^k}{c-q^k}q^k\bigg)\\\times
\frac{\Big(\frac{1-bq^k}{c-q^k}dq;q\Big)_\infty
\Big(\frac{c-q^k}{1-bq^k}eq;q\Big)_\infty}
{\Big(\frac{1-bq^k}{c-q^k};q\Big)_\infty
\Big(\frac{c-q^k}{1-bq^k}q;q\Big)_\infty}\,
\frac{(1/d;q)_k}{(eq;q)_k}
\left(\frac{1-bq^k}{c-q^k}d\right)^k,
\end{multline}
where $|d/c|<1$ and $|e/b|<1$.
\end{Theorem}

The proof is similar to the proof of Theorem~\ref{bns}.
We therefore omit the details.
Theorem~\ref{bnsc} reduces to the $a=0$ case of
Theorem~\ref{tnsc} for $e=1$. Similar to \eqref{bnsgl},
the $c=1/b$ case of \eqref{bnscgl} gives,
after some rewriting, Ramanujan's ${}_1\psi_1$ summation \eqref{11s}.
Also, after multiplying both sides by $c$,
performing the substitution $c\mapsto -b/c$ and
then letting $b\to\infty$, one obtains the very-well-poised
${}_4\psi_6$ summation in \eqref{46s}.

\end{document}